\documentclass[12pt]{article}
\usepackage[english]{babel}
\usepackage[centertags]{amsmath}
\usepackage{amsfonts}
\usepackage{amssymb}
\usepackage{amsthm}
\usepackage{newlfont}
\usepackage[bookmarksnumbered, colorlinks, plainpages]{hyperref}

\usepackage{natbib}
\begin{document}

\title{Some estimates for the norm of the  self-commutator}



\author{Gevorgyan Levon}



\maketitle
{Key words: Self-commutator, Numerical range}

AMS Mathematics Subject Classification {47A12,47B47 \{Primary\}, 52A10, 52A40 \{Secondary\}}

\begin{abstract}
Different estimates for the norm of the self-commutator of a Hilbert space operator are proposed. Particularly, this norm is bounded from above by twice of the area of the numerical range of the operator. An isoperimetric-type inequality is proved. 

\end{abstract}
\maketitle
\section{Introduction and preliminaries}

\noindent Let  $A$ be a linear bounded operator, acting in a Hilbert
space $\left ( {\mathcal{H},\left\langle { \bullet , \bullet }
\right\rangle }\right).$ The difference $ A^*A  - AA^* = C\left( A \right)$ is said to
be the self-commutator of the operator $A.$ The study of self-commutators was initiated by  Halmos in \citet{MR0045310}. The well-known and important class of normal operators is characterized by the equality $ AA^*  = A^*A,$    so the norm of $C\left( A \right)$ shows "how far" is the operator from being normal.
If $C\left(A \right)$ is semi-definite, the operator $A$ is said \citet{MR0217618} to be
semi-normal, particularly, if $C\left( A \right) \ge \textbf{0},$
then $A$ is hyponormal. According to Putnam's inequality
for any semi-normal operator
\[\left\| {C\left( A\right)} \right\| \le \frac{1}{\pi }mes_2 \left( {SpA} \right),\]
where $SpA$  is the spectrum of $A$  and $mes_2$ means the plane Lebesgue measure.

In \citet{MR2655007} is proved that for any operator $A$ the inequality \begin{equation}
\left\| C\left( A\right) \right\| \le \left\| A \right\|^2.\label{one}
\end{equation}
is satisfied and  an operator $A$ is said  to have large self-commutator if  \newline
$\left\| C\left( A\right)\right\| = \left\| A \right\|^2 .$ As $C\left( {A - \lambda I} \right) = C\left( A\right)$ for any $\lambda  \in \mathbb{C},$ inequality \eqref{one} may be sharpened \[\left\| C\left( A\right) \right\| \leq m^2 \left( A \right),\]
where \[m\left( A \right) = \mathop {\inf }\limits_{\lambda  \in \mathbb{C}} \left\| {A - \lambda I} \right\|.\]

In the sequel we propose  less obvious estimates for the norm of the self-commutator of an operator. Particularly, $C\left(A\right)$  will be bounded from above by the twice of the area of the numerical range of $W\left(A\right)$ of $A.$ For operators with elliptical numerical range the factor 2 may be replced by $4/\pi.$

\section{Main results}
 \textbf{\emph{Proposition 1.}} Let $N\left( A \right)$ the set of all normal operators, commuting with $A.$ Then \begin{equation}
\left\| C\left( A\right) \right\| \leqslant \mathop {\inf }\limits_{M \in N\left( A \right)} \left\| {A + M} \right\|^2. \label{PF}\end{equation}

\textbf{\emph{Proof.}} Let $M \in N\left( A \right).$  Then it is easy to see that $C\left( A \right) = C\left( {A + M} \right).$ Indeed, according to the Putnam-Fuglede theorem (\citet{MR675952}, Problem 152) the equality $AM = MA$ implies
$A^* M = MA^* , AM^*  = M^* A,$ hence \newline
$\left( {A^*  + M^* } \right)\left( {A + M} \right) - \left( {A + M} \right)\left( {A^*  + M^* } \right) = A^* A - AA^*.$ The proof may be completed recalling formula \eqref{one}.

\textbf{\emph{Remark.}} From this Proposition follows that for any hyponormal operator $A$ the operator$A+M, M \in N\left( A \right)$ is also hyponormal, so in the Putnam  inequalities at the right-hand side one may put the area of the spectrum of any operator $A+M,$ which may have smaller spectrum.

An operator $A$ is said to be irreducible, if the only invariant under $A$ subspaces are two trivial subspaces $\left\{ \theta  \right\}$ and $ \mathcal{H},$ or which is the same, the only orthogonal projections, commuting with $A$ are $O$ and $I.$

\emph{\textbf{Example 1.}} The operator of the simple unilateral shift $U$ is irreducible (\citet{MR675952}, Problem 116, Corollary),

\textbf{\emph{Proposition 2.}} Let $A$ be irreducible. Then \begin{equation}
\mathop {\inf }\limits_{M \in N\left( A \right)} \left\| {A + M} \right\| = \mathop {\inf }\limits_{\lambda  \in \mathbb{C}} \left\| {A - \lambda I} \right\|.\label{occ}\end{equation}

\textbf{\emph{Proof.}} Let $M \in N\left( A \right).$  As the spectral projections $\left\{ P \right\}$ of $M$ commute with any operator, commuting with $M,$ the spectral projection $P$ commutes with $A,$ hence $P = 0$ or $P=I,$ implying $M=\lambda I,\,\,\lambda  \in \mathbb{C}.$

The next example shows that the left hand side in  \eqref{occ} may be strictly less than other side.

\emph{\textbf{Example 2.}} Let
\[L = \left( {\begin{array}{*{20}c}
   1 & 1 & 0 & 0  \\
   0 & 1 & 0 & 0  \\
   0 & 0 & { - 1} & 1  \\
   0 & 0 & 0 & { - 1}  \\

 \end{array} } \right)\]
and
 \[K = \left( {\begin{array}{*{20}c}
   1 & 0 & 0 & 0  \\
   0 & 1 & 0 & 0  \\
   0 & 0 & { - 1} & 0  \\
   0 & 0 & 0 & { - 1}  \\

 \end{array} } \right).\]
Evidently
\[KL = LK,\,K^*  = K,\left\| {L - K} \right\| = 1,\left\| {L - \lambda I} \right\| \geqslant \frac{{1 + \sqrt 5 }}
{2}.\]

For the norm of $C\left( A \right)$ there are some estimates from below.
In finite dimensional space  it is shown (\citet{MR822142}, Proposition 2) that there exists a matrix $B$ such that
$C\left( B \right) = C\left( A \right)$ and $\left\| {C\left( A \right)} \right\| \geqslant \left\| B \right\|^2 /2.$

The matrix $B$ is constructed by the following procedure.
Let $A$ be  \[
A = \left( {\begin{array}{*{20}c}
   0 & 0 & 0 &  \cdots  & 0  \\
   {\xi _1 } & 0 & 0 &  \cdots  & 0  \\
   0 & {\xi _2 } & 0 &  \cdots  & 0  \\
    \cdot  &  \cdot  &  \cdot  &  \cdot  &  \cdot   \\
   0 & 0 &  \cdots  & {\xi _{n - 1} } & 0  \\
 \end{array} } \right).\]
It is easy to see that $\left\| A \right\| = \mathop {\max }\limits_k \left| {\xi _k } \right|.$ The operator $C\left(A \right)$ has the following form \[
diag\left( {\left| {\xi _1 } \right|^2 ,\left| {\xi _2 } \right|^2  - \left| {\xi _1 } \right|^2 , \cdots ,\left| {\xi _{k + 1} } \right|^2  - \left| {\xi _k } \right|^2 , \cdots , - \left| {\xi _{n - 1} } \right|^2 } \right).\]

Therefore the eigenvalues $\left\{ {\lambda _k } \right\}_1^n$ of $C\left(A \right)$ satisfy the relations \begin{equation}
\sum\limits_{k = 1}^m {\lambda _k }  = \left| {\xi _m } \right|^2  \geqslant 0,\,\left( {1 \leqslant m \leqslant n - 1} \right).
\label{Fo}\end{equation}
We have \[
\left\| {C\left( A \right)} \right\| = \mathop {\max }\limits_k \left| {\lambda _k } \right|.
\]
For any set of real numbers $\left\{ {\lambda _k } \right\}_1^n$ having zero sum it is possible to rearrange them such that condition \eqref{Fo} is satisfied and \begin{equation}
\sum\limits_{k = 1}^m {\lambda _k }  \leqslant 2\mathop {\max }\limits_k \left| {\lambda _k } \right|.
\label{bin}\end{equation}
Finally we get \[
\left\| {C\left( A \right)} \right\| \geqslant \left\| A \right\|^2 /2.\]

The same may be repeated in some infinite dimensional spaces.

Let $A^2\left(D\right)$ be the Bergman space of analytic in the unit circle functions square integrable with respect to the plane Lebesgue measure in $D.$ The set of monomials $\left\{ {e_n=\sqrt {n + 1} z^n } \right\}_0^\infty$ is an orthonormal basis of $A^2\left(D\right).$ Denote by $M$ the operator of multiplication by the independent variable \[\left( {Mf} \right)\left( z \right) = zf\left( z \right).\]
We have
\[\left( M^{*}f\right) \left( z\right) =z^{-2}\left( zf\left( z\right)
-\int\limits_0^zf\left( \xi \right) d\xi \right).\]
It is easy to see that $M$ is the operator of the weighted shift $Me_n=w_ne_{n+1}$ with the weight sequence $w_n  = \sqrt {\frac{{n + 1}}
{{n + 2}}}.$
The operator $C\left(M \right)$ is defined by the formula
\[\left( {\left( {M^* M - MM^* } \right)f} \right)\left( z \right) = \frac{1}
{{z^2 }}\int\limits_0^z {\left( {z - \xi } \right)f\left( \xi  \right)d\xi }.\]
The monomials $\left\{ {z^n } \right\}_0^\infty$ are the eigenfunctions of this operator, corresponding to the eigenvalues
$\left\{ {\frac{1}
{{\left( {n + 1} \right)\left( {n + 2} \right)}}} \right\}_0^\infty.$
One has $\left\|M\right\|=1, \left\|C\left(M\right)\right\|=1/2.$

Note that
\[0 \le \sum\limits_{k = 0}^m {\frac{1}
{{\left( {k + 1} \right)\left( {k + 2} \right)}}}  < 1 = 2\sup \lambda _n ,\]
hence inequalities \eqref{Fo}, \eqref{bin} for this operator are satisfied.

Now we pass to the general case.
Denote $A_t  = Ae^{it} ,t \in \left[ {0;2\pi } \right)$ and
$H_t  = \operatorname{Re} A_t, J_t  = \operatorname{Im} A_t , H=H_0, J=J_0.$ Note that $H_{t + \pi /2}  =  - J_t,$ implying
\[ \mathop {\max }\limits_{t \in \left[ {0;2\pi } \right)} \left\| {J_t } \right\| = \mathop {\max }\limits_{t \in \left[ {0;2\pi } \right)} \left\| {H_t } \right\|=w\left( A \right),\]
where  $w\left( A \right)$ is the numerical radius of the operator  $A.$

Denote \[\alpha _2  = \mathop {\sup }\limits_{\left\| x \right\| = 1} \left\langle {H_t x,x} \right\rangle ,\alpha _1  = \mathop {\inf }\limits_{\left\| x \right\| = 1} \left\langle {H_t x,x} \right\rangle ,\]
\[\beta _2  = \mathop {\sup }\limits_{\left\| x \right\| = 1} \left\langle {J_t x,x} \right\rangle ,\beta _1  = \mathop {\inf }\limits_{\left\| x \right\| = 1} \left\langle {J_t x,x} \right\rangle .\]
It fact $b_x \left(t \right) = \alpha_2-\alpha_1$ ($b_y\left(t \right) =\beta_2-\beta_1$) is the width of $W \left( A_t \right)$ in the direction of $Ox$-axis ($Oy$-axis).

\textbf{\emph{Remark.}} As the closure of the numerical range of a self-adjoint operator coincides with the convex hull of its spectrum, the infimum and the supremum over the numerical ranges may be replaced by corresponding bounds over the spectrum of respective operators.

\textbf{\emph{Proposition 3.}} For any Hilbert space operator $A$ the following inequality is satisfied

\begin{equation}\left\| {C\left( A \right)} \right\| \leqslant \mathop {\min }\limits_{t \in \left[ {0;2\pi } \right)} \left\{ {b_x \left( t \right)b_y \left( t \right)} \right\}.\label{my}\end{equation}

\textbf{\emph{Proof.}} Easy to check that the self-commutator is rotation-invariant, hence
 $\left\| {C\left( A \right)} \right\| = 2\left\| {J_t H_t  - H_t J} \right\|.$

According to \citet{MR2440836}, Corollary 7
\[\left\| {J_t H_t  - H_t J_t } \right\| \leqslant \frac{1}
{2}b_x \left( t \right)b_y \left( t \right)\]
and  from arbitrariness of $t$ follows \eqref{my}.

In \citet{MR2440836}, Corollary 11  the inequality
\begin{equation}\left\| {A^* A - AA^* } \right\| \leqslant 4\mathop {\inf }\limits_{\lambda  \in \mathbb{C}} \left\| {H - \lambda I} \right\|\cdot\mathop {\inf }\limits_{\lambda  \in \mathbb{C}} \left\| {J - \lambda I} \right\|\label{WaDu}\end{equation}
is proved.

The following example shows that \eqref{my} is really sharper than \eqref{WaDu}.

\textbf{\emph{Example 3.}} Let
 \[A = \frac{{\sqrt 2 }}
{2}\left( {\begin{array}{*{20}c}
   {1 + i} & {1 + i}  \\
   0 & { - 1 - i}  \\
 \end{array} } \right).\]
Then
\[H = J^* = \frac{{\sqrt 2 }}
{2}\left( {\begin{array}{*{20}c}
   1 & {\frac{{1 + i}}
{2}}  \\
   {\frac{{1 - i}}
{2}} & { - 1}  \\

 \end{array} } \right).\]
We have $\left\| A \right\| = \sqrt 5, $
\[\mathop {\inf }\limits_{\lambda  \in \mathbb{C}} \left\| {H - \lambda I} \right\| = \mathop {\inf }\limits_{\lambda  \in \mathbb{C}} \left\| {J - \lambda I} \right\| = \frac{{\sqrt 3 }}{2}\]
and $\mathop {\min }\limits_{t \in \left[ {0;2\pi } \right)} \left\{ {b_x \left( t \right)b_y \left( t \right)} \right\}=\sqrt 5.$

Combining Proposition 1 from the Appendix and Proposition 3, we get the following result.

\textbf{\emph{Proposition 4.}} For any operator $A$
\begin{equation}
\left\| {C\left( A \right)} \right\| \leqslant 2S\left( {W\left( A \right)} \right).
\label{alm}\end{equation}
 For some operators estimate \eqref{alm} may be sharpened.

\textbf{\emph{Proposition 5.}} Let $A$ have elliptical numerical range. Then
\begin{equation}
\left\| {C\left( A \right)} \right\| \leqslant \frac{4}
{\pi }S\left( W \left( A\right)\right).\label{censym}
\end{equation}

\textbf{\emph{Proof.}} As both the self-commutator and the area of the numerical range are translation invariant, we may suppose that the centre of the ellipse coincides with the origin of the coordinate system. Recall now Proposition 3 and Example 2 (formula \eqref{ellips}) from the Appendix.

\section{Examples}
\textbf{\emph{Example 1.}} Let the operator $A$ be the tensor product of $a,b \in H,$ i.e $Ax = \left\langle {x,a} \right\rangle b.$ To avoid triviality, we exclude the case of $a = \theta$ and $b=\theta.$
We intend to calculate the norm of the self-commutator of this operator. If $a$ and $b$ are parallel, the $A$ is normal, hence its self-commutator is equal to zero. This possibility will be also ruled out.  Evidently $\left\| A \right\| = \left\| a \right\| \cdot \left\| b \right\|.$  As $
C\left( {\lambda A} \right) = \left| \lambda  \right|^2 C\left( A \right)$ then without loss of generality  we may assume that  $\left\| A \right\| = 1,$  hence $\left\| a \right\| = \left\| b \right\| = 1.$ Then $A^* Ax = \left\langle {x,a} \right\rangle a = P_a x,\,AA^* x = \left\langle {x,b} \right\rangle b = P_b x,$
where $P_ \bullet$ is the operator of the orthogonal projection onto the subspace, generated  by $\bullet.$ Thus \[
\left( {A^* A - AA^* } \right)x = \left( {P_a  - P_b } \right)x.\]  Recall that the norm of $\left\| {P_a  - P_b } \right\|$ is said to be  the opening of two subspaces, generated by $a$ and $b$ respectively. Denote by  $\bigvee ab$ the linear span of $a$ and $b.$
We get the following result.

\textbf{\emph{Proposition 1.}} \begin{equation}
\left\| {P_a  - P_b } \right\| = \left\| x \right\|\sqrt {1 - \left| {\left\langle {a,b} \right\rangle } \right|^2.}\label{Mit}\end{equation}

\textbf{\emph{Proof.}} Let  $\left\{ {a_k } \right\}_1^n$ be a set of linear independent elements from $H$ and $G$ be their Gram matrix. It is known (the Generalized Ostrowski's inequality, \citet{MR2085667})  that \citet{MR2302906}, ch. III, 93
\begin{equation}\left\| x \right\|^2  \ge \sum\limits_{i,k = 1}^n {g_{ik}^{ - 1} } \left\langle {x,a_i } \right\rangle \left\langle {a_k ,x} \right\rangle, \label{glaz}\end{equation}
where $g_{ik}^{ - 1}$ are the elements of the matrix $G^{-1}.$  The equality is attained if and only if  $x$ belongs to the linear span of $\left\{ {a_k } \right\}_1^n.$

For $n=2$ we get
\[G = \left( {\begin{array}{*{20}c}
   1 & {\left\langle {a,b} \right\rangle }  \\
   {\left\langle {b,a} \right\rangle } & 1  \\
\end{array}} \right)\]
and
\[G^{ - 1}  = \frac{1}
{{1 - \left| {\left\langle {a,b} \right\rangle } \right|^2 }}\left( {\begin{array}{*{20}c}
   1 & { - \left\langle {a,b} \right\rangle }  \\
   { - \left\langle {b,a} \right\rangle } & 1  \\

 \end{array} } \right)\]
implying
\begin{equation}
\begin{gathered}
  \left\| x \right\|^2  \geqslant \left\| {P_{ \bigvee ab} x} \right\|^2  = \frac{1}
{{1 - \left| {\left\langle {a,b} \right\rangle } \right|^2 }} \times  \hfill \\
   \times \left( {\left| {\left\langle {x,a} \right\rangle } \right|^2  - \left\langle {a,b} \right\rangle \left\langle {x,a} \right\rangle \left\langle {b,x} \right\rangle  - \left\langle {b,a} \right\rangle \left\langle {x,b} \right\rangle \left\langle {a,x} \right\rangle  + \left| {\left\langle {x,b} \right\rangle } \right|^2 } \right). \hfill \\
\end{gathered}\label{huge}\end{equation}

As the expression in the parentheses is equal to $\left\| {\left\langle {x,a} \right\rangle a - \left\langle {x,b} \right\rangle b} \right\|^2$ we have\[
\left\| {P_a x - P_b x} \right\| \leqslant \left\| x \right\| \cdot \sqrt {1 - \left| {\left\langle {a,b} \right\rangle } \right|^2 },
\]
implying \eqref{Mit}.

Dragomir \citet{MR2366098} proved that for any operator $A$ the following inequality \[
\left\| A \right\|^2  - w^2 \left( A \right) \leqslant m^2 \left( A \right)\]
is satisfied. It is easy to see that for considered above operator a more accurate inequality \begin{equation}
\left\| A \right\|^2  - w^2 \left( A \right) \leqslant C\left( A \right)\label{hypo}\end{equation} takes place.

\textbf{\emph{Conjecture 1.}} For any Hilbert space operator $A$ inequality \eqref{hypo} is satisfied.

\textbf{\emph{Example 2.}} Let $A$ be arbitrary operator, acting in a two dimensional space. It may be reduced to the Schur's upper triangular form \[
A = \left( {\begin{array}{*{20}c}
   {\lambda _1 } & {\lambda _3 }  \\
   0 & {\lambda _2 }  \\

 \end{array} } \right).\]

Recall that the numerical range of a $2\times2$ matrix is an elliptical disk (\citet{MR675952}, ch. 22) having as foci two
eigenvalues $\left\{ {\lambda _1 ,\lambda _2 } \right\}$ and \[
\left\{ {tr\left( {A^* A} \right) - \left| {\lambda _1 } \right|^2  - \left| {\lambda _2 } \right|^2 } \right\}^{1/2}\] as minor axis $2b.$

We have  $b = \left| {\lambda _3 } \right|/2,$ the distance between foci is
$2c = \left| {\lambda _2  - \lambda _1 } \right|,$ hence major axis is
$2a = \sqrt {\left| {\lambda _2  - \lambda _1 } \right|^2  + \left| {\lambda _3 } \right|^2 }.$
 Easy calculations show that \[
A^* A - AA^*  = \left( {\begin{array}{*{20}c}
   { - \left| {\lambda _3 } \right|^2 } & {\lambda _3 \left( {\overline \lambda  _1  - \overline \lambda  _2 } \right)}  \\
   {\overline \lambda  _3 \left( {\lambda _1  - \lambda _2 } \right)} & {\left| {\lambda _3 } \right|^2 }  \\

 \end{array} } \right).\]
The norm of this operator is equal to $\left| {\lambda _3 } \right|\sqrt {\left| {\lambda _2  - \lambda _1 } \right|^2  + \left| {\lambda _3 } \right|^2}$ or   $4ab = 4S/\pi,$ where $S$ is the area of the ellipse-the numerical range of $A.$

\textbf{\emph{Example 3.}} Let $A$ be the SOR iteration matrix
\[A = \left( {\begin{array}{*{20}c}
   {\left( {1 - \omega } \right)I_p } & {\omega M}  \\
   {\omega \left( {1 - \omega } \right)M^T } & {\left( {1 - \omega } \right)I_q  + \omega ^2 M^T M}  \\

 \end{array} } \right),\]
where $0 < \omega  < 2,\;\;M \in \mathbb{R}^{p \times q} ,\;p \geqslant q$ and $I$ is the identity matrix. Golub and de Pillis showed in \citet{MR1038083} that $A$ is unitary equivalent to the matrix \newline $diag\left\{ {M_1 ,M_2 , \cdots ,M_q ,\left( {1 - \omega } \right)I_{p - q} } \right\},$  where $M_k  \in \mathbb{R}^{2 \times 2}.$

Evidently $C\left(A\right)= diag\left\{ {C\left( {M_1 } \right),C\left( {M_2 } \right), \cdots ,C\left( {M_q } \right),0} \right\}$ and $\left\| {C\left( A \right)} \right\| = \max \left\| {C\left( {M_k } \right)} \right\|.$ As $
W\left( A \right) = ch\left\{ {W\left( {M_1 } \right),W\left( {M_2 } \right), \cdots ,W\left( {M_q } \right),1 - \omega } \right\}$ and $
1 - \omega  \in W\left( {M_k } \right),1 \leqslant k \leqslant q,$ we get
\[
\left\| {C\left( A \right)} \right\| \leqslant \frac{4}
{\pi }S\left( {W\left( A \right)} \right).
\]

\textbf{\emph{Example 4.}} Let $\emph{D}$ be the Dirichlet space of the functions, analytic in the open unit circle and having finite norm
\[\left\| f \right\|^2  = \left| {f(0)} \right|^2  + \frac{1}
{\pi }\iint\limits_{x^2  + y^2  < 1} {\left| {f'} \right|}^2 dxdy\]
and
$\phi  = \phi \left( {z,a} \right) = \frac{{a-z}}{{1 - \overline a z}},\left| a \right| < 1$
be the M\"{o}bius function. Denote by $C_\phi$ the composition operator, induced by $\phi,$ i.e. $C_\phi  f = f \circ \phi.$ It is known \citet{MR2407083} that $\left\| {C_\phi ^* C_\phi   - C_\phi  C_\phi ^* } \right\| = \sqrt {L^2  + 4L} ,$ where $L =  - \ln \left( {1 - \left| a \right|^2 } \right).$ The closure of the numerical range of $C_\phi $ is an ellipse with foci $\pm 1$ and with the major axis length $\sqrt {4 + L}.$ The minor axis length is $2b =  \sqrt L,$ therefore the area of the ellipse is
$S = \frac{\pi }
{4}2a2b = \frac{\pi }
{4}\sqrt {L^2  + 4L}$
and
$\left\| {C_\phi ^* C_\phi   - C_\phi  C_\phi ^* } \right\| = 4S/\pi.$

\textbf{\emph{Example 5.}} Let $A$ be a tridiagonal Toeplitz matrix
\[A = \left( {\begin{array}{*{20}c}
   \lambda  & a & 0 &  \cdots  & 0 & 0  \\
   b & \lambda  & a &  \cdots  & 0 & 0  \\
    \cdot  &  \cdot  &  \cdot  &  \cdot  &  \cdot  &  \cdot   \\
   0 & 0 & 0 &  \cdots  & \lambda  & a  \\
   0 & 0 & 0 &  \cdots  & b & \lambda   \\

 \end{array} } \right),\]
where $a,b,\lambda  \in \mathbb{C}.$ It is easy to see that for the self-commutator
$C = C\left( A \right)$ all the elements are equal to zero, except $
C_{11}  =  - C_{nn}  = \left| b \right|^2  - \left| a \right|^2 ,
$ implying $\left\| C \right\| = \left| {\left| a \right|^2  - \left| b \right|^2 } \right|.$

The numerical range $W\left( {A - \lambda I} \right)$ is \citet{MR1206415} the ellipse, defined by the formula
$z = \left( {ae^{ - it}  + be^{it} } \right)\cos \frac{\pi }
{{n + 1}},t \in \left[ {0,2\pi } \right].$
As the semi-axes of this ellipse are
$\left( {\left| a \right| + \left| b \right|} \right)\cos \frac{\pi }
{{n + 1}},\,\left| {\left| a \right| - \left| b \right|} \right|\cos \frac{\pi }
{{n + 1}},$
the area $S$ of the ellipse is $S = \pi \,\left| {\left| a \right|^2  - \left| b \right|^2 } \right|\cos ^2 \frac{\pi }
{{n + 1}}.$ Using the inequality
$\cos ^2 \frac{\pi }
{{n + 1}} \geqslant \frac{1}{4}$
for $n \geqslant 2,$ we get finally \eqref{censym}.

\textbf{\emph{Example 6.}} Consider the Volterra integration operator in $
L^2 \left( {0\,;1} \right) $ defined by the formula

\[\left( {Vf} \right)\left( x \right) = \int\limits_0^x
{f\left( t \right)dt}.\]
Recall (\citet{MR675952}, Problem 165) that $ W\left( V \right)$ is bounded by the curve
\[t \mapsto \frac{{1 - \cos t}}{{t^2 }} \pm i\frac{{t - \sin
t}}{{t^2 }},\,0 \le t \le 2\pi.\]
The area of $W\left(V\right)$ may be calculated by the formula
\[S = \int\limits_0^{2\pi } {xdy - ydx}  = \frac{1}
{6}Si\left( {2\pi } \right) + \frac{1}
{{12\pi }} \approx 0.26288441987 \cdots ,\]
where $Si$ is the integral sine function
\[Si\left( x \right) = \int\limits_0^x {\frac{{\sin t}}
{t}dt}.\]
We get $\frac{4}
{\pi }S \approx 0.33471483907 \cdots$
On the other hand $\left\| {C\left( V \right)} \right\|{\text{ }} = \frac{{\sqrt 3 }}
{6} \approx 0.288675134...$

\textbf{\emph{Conjecture 2.}} For any Hilbert space operator $A$ inequality \eqref{censym} is satisfied.

\section{Appendix}

Let $G$ be a convex compact figure in the plane $\mathbb{C}$ and $S\left(G\right)-$ the area of $G.$ Denote $
k = \left\{ {\cos t;\sin t} \right\},j = \left\{ { - \sin t;\cos t} \right\},t \in \left[ {0;2\pi } \right)$
 and let $w_k \left( G \right)$,
 \newline$\left( {w_j \left( G \right)} \right)$ be the width of $G$ in the direction of $k$ ($j$). Denote by $d$ the diameter of $G, l\left( d \right)-$ the length of $d.$  According to (\citet{MR0221383}, 11, Properties 1-6) there exist two line segments $d$ and $h$ connecting points of the boundary of $G,$ the directions of which are the directions of greatest and least width, perpendicular to the supporting lines drawn through their endpoints.

\textbf{\emph{Proposition 1.}} For any $G$ the following inequality is satisfied
\begin{equation}
\mathop {\inf }\limits_{t \in \left[ {0;2\pi } \right)} \left\{ {w_j \left( G \right)w_k \left( G \right)} \right\} \leqslant 2S\left( G \right).
\label{inv}\end{equation}

\textbf{\emph{Proof.}} We prove that the inequality above is satisfied for a particular value of $t.$ Let $t_0$ be chosen in a such way that $j$ has the direction of the least width $h.$ Denote by $n$ and $s$ two extremities of $h$ and by $e$ and $w-$ farthest from $h$ points of $G$ lying in the different half-planes, separated by $h.$ The quadrilateral $Q$ with vertices $nesw$ is a subset of $G$ hence $S\left( Q \right) \leqslant S\left( G \right).$
On the other hand $w_j \left( G \right)w_k \left( G \right) = w_j \left( Q \right)w_k \left( Q \right) = 2S\left( Q \right),$  completing the proof.

For any convex set $G$ \begin{equation}
1 \leqslant \frac{1}
{{S\left( G \right)}}\mathop {\inf }\limits_{t \in \left[ {0;2\pi } \right)} \left\{ {w_j \left( G \right)w_k \left( G \right)} \right\} \leqslant 2.\label{ratio}
\end{equation}
The lower bound is attained for a rectangle. The following example shows that the upper bound is also exact.

\emph{\textbf{Example 1.}} Consider the equilateral triangle with the summits $z_k  = \exp \left( {i\frac{{2\pi k}}{3}} \right),k = 0,1,2.$ By symmetry it is sufficient to consider the case $0 \leqslant t \leqslant \frac{\pi }
{6}.$  We have $w_k w_j  = 3\sin \left( {\frac{\pi }
{3} + t} \right)\cos t.$ The minimum of this product is $\frac{{3\sqrt 3 }}{2}$ and is attained at the endpoints of mentioned above segment.

For particular convex sets the constant $2$ may be diminished.

\emph{\textbf{Example 2.}} Let $G$ be the ellipse, defined by the canonical equation $\frac{{x^2 }}
{{a^2 }} + \frac{{y^2 }}
{{b^2 }} = 1.$
Then (\citet{MR1794773}, Proposition 2.3) $w_k \left( t \right) = 2\sqrt {a^2 \cos ^2 t + b\sin ^2 t}$,$w_j \left( t \right) = 2\sqrt {a^2 \sin ^2 t + b^2 \cos ^2 t}.$ Thus \begin{equation}
\mathop {\min }\limits_{t \in \left[ {0;2\pi } \right)} \left\{ {w_k \left( G \right)w_j \left( G \right)} \right\} = 4ab = \frac{4}
{\pi }S\left( G \right).\label{ellips}\end{equation}

This estimate holds also for the particular case - a disk. For any figure of the constant width the  ratio \eqref{ratio} (according to Blaschke-Lebesgue theorem, \citet{gruber1983convexity}) lies between two extremal values $4/\pi\approx 1.273239545 \cdots$ and $\frac{2}
{{\pi  - \sqrt 3 }} \approx 1.418900762 \cdots,$ which is the ratio  for the Reuleaux triangle.

\bibliographystyle{amsplain}

\bibliography{mybib1}
\end{document}